\documentclass{amsart}
\usepackage{amsmath}
  \usepackage{paralist}
 \usepackage{hyperref}
   
 \usepackage{tikz}
\usepackage{graphicx}
\usepackage{amssymb}
\usepackage{epstopdf}
\usepackage{amsthm}
\usepackage{txfonts}
\usepackage{mathrsfs}

\usepackage[round]{natbib}

\usepackage{xfrac}
\usepackage{multicol}
\usepackage{mathabx}
\usepackage{tabularx}
\usepackage{array}
\usepackage{flafter}
\usepackage{float}
\usepackage{subcaption}
\usepackage{wrapfig}
\usepackage[symbol*]{footmisc}
\usepackage{arydshln}
\DefineFNsymbolsTM{myfnsymbols}{
  \textasteriskcentered *
  \textdagger    \dagger
  \textdaggerdbl \ddagger
  \textsection   \mathsection
  \textbardbl    \|%
  \textparagraph \mathparagraph
}%
\setfnsymbol{myfnsymbols}

\floatstyle{boxed} 
\restylefloat{figure}

\newtheorem{thm}{Theorem}[section]

\theoremstyle{definition}

\newtheorem{rmk}{Remark}

\newcommand{\p}{\\[10pt]\noindent}

\renewcommand{\bar}{\overline}

\newcommand{\al}{\alpha}
\newcommand{\om}{\omega}

\newcommand{\hlf}{\frac{1}{2}}

\renewcommand{\phi}{\varphi}

\newcommand{\mR}{\mathcal{R}}

\newcommand*{\myproofname}{Proof}

\newcommand{\LeftEqNo}{\let\veqno\@@leqno}

\title[Extinction Probabilities for ISAv] 
      {The Probability of Extinction of Infectious Salmon Anemia Virus in One and Two Patches}

\author[Evan Milliken]{}

\thanks{The author is supported by NSF grants DMS-1411853, DMS-1515661 and the Center for Applied Mathematics at University of Florida}

\begin{document}
\maketitle

\centerline{\scshape Evan Milliken$^*$}
\medskip
{\footnotesize
 \centerline{School of Mathematical and Statistical Sciences, Arizona State University}

} 

%

\bigskip


\begin{abstract}  Single type and multitype branching process have been used to study the dynamics of a variety of stochastic birth-death type phenomena in biology and physics.  Their use in epidemiology goes back to Whittle's study of a Susceptible--Infected--Recovered (SIR) model in the 1950s.  In the case of an SIR model, the presence of only one infectious class allows for the use of single type branching processes.  Multitype branching processes allow for multiple infectious classes and have latterly been used to study metapopulation models of disease.  In this article, we develop a Continuous Time Markov Chain (CTMC) model of Infectious Salmon Anemia virus in two patches, two CTMC models in one patch and companion multitype branching process (MTBP) models.  The CTMC models are related to deterministic models which inform the choice of parameters.  The probability of extinction is computed for the CTMC via numerical methods and approximated by the MTBP in the supercritical regime.  The stochastic models are treated as toy models and the parameter choices are made to highlight regions of the parameter space where CTMC and MTBP agree or disagree, without regard to biological significance.  Partial extinction events are defined and their relevance discussed.  A case is made for calculating the probability of such events, noting that MTBP's are not suitable for making these calculations.
\end{abstract}

\section{Introduction}
In the investigation that follows, we will use the case of an outbreak of Infectious Salmon Anemia (ISA) as a test case to examine some of the features of MTBP approximation of a CTMC.  Infectious salmon anemia virus (ISAv) causes (ISA) which leads to 15 -- 100\% accumulated mortality over the course of a several months long infection in a farm environment \cite{Falk1997}.  It is found in all large salmon-producing countries including Norway, Scotland, Ireland, Canada, the United States, and Chile \cite{Vike2009}.  ISAv is transmitted among finfish horizontally by passive movement of infected seawater \cite{Mardones2009} and via direct contact with excretions or secretions of infected individuals.  Salmon farms consist of a collection of net cages placed in open body of water.  This array-like structure of a farm and the proximity of farms to each other and to wild salmon migratory routes justifies the use of a metapopulation approach.

Branching processes have been used to study a variety of biological phenomena dating back to their invention to answer a question regarding the extinction of aristocratic surnames.   Bienaym\'{e} made the first contribution in 1845 \cite{Seneta1998} before the question was made well known by Galton and answered together with Watson in 1873-4 \cite{Watson1875}.  As a result, the class of single type branching processes came to be known as Bienaym\'{e}--Galton--Watson branching processes (BGWbp).  A special case considering two types of individuals was studied by Bartlett in 1946 and BGWbp theory was extended to include general multitype branching processes by Kolmogorov, Dmitriev, Sevastyanov, Everett and Ulam in the late 1940's \cite{Harris1963}.  BGWbp and MTBP models have been used to study a variety of phenomena in biology and physics including population dynamics, changes to the genome, cell kinetics, cancer and epidemiology \cite{Allen2003,Allen2012,Allen2013a,Allen2013,Ball1983,Ball1995,Britton2010,Dorman2004,Griffiths2011,Harris1963,Kimmel2002,Whittle1955}.  In particular, Allen and Lahodny studied MTBP's as an approximation of the outbreak dynamics of a CTMC model of infection in single and multi-patch models \cite{Allen2012,Allen2013a}.

We recall earlier analysis of deterministic Susceptible-Infected-Virus (SIV) models of ISAv outbreak in one and two patches to inform our investigation and aid in suitable parameter selection \cite{Milliken2016}.  The model in one patch is adapted from well-studied models \cite{Beretta1998,Nowak2000,Perelson1999} by allowing for direct transmission via contact with infected individuals.  For each of these two models, a companion CTMC model is introduced, as well as a MTBP.  The probability of extinction of the disease is approximated for the CTMC using numerical simulation.  Approximation is also made via analysis of the MTBP and the results are compared with those of numerical simulation.

Formulation of a birth--death process as a branching process relies on the fact that all transitions are independent \cite{Allen2012,Harris1963,Mode1971}.  This is a strong biological assumption, but one commonly made for the purpose of mathematical modeling.  In order to formulate epidemiological models as branching processes, an additional assumption is made:  the susceptible population remains fixed at its initial (disease-free) population size.  As a result of this assumption, MTBP only provides accurate approximation of the probability of disease extinction when the total population size is sufficiently large.  There is currently no analytic estimate for how large is sufficiently large.  In order to illustrate the breakdown of MTBP approximation and explore its dependence on the underlying system and its parameters, we calculate the probability of extinction for a one-patch system for a range of initial population sizes at two different levels of infected fish mortality.  We also propose a variation on the deterministic one-patch model by changing the assumed force of infection (\emph{f.o.i.}).  Corresponding CTMC and MTBP models are also developed.  The probability of extinction is again calculated at various initial population sizes.

MTBP techniques are suitable to calculate the probability of complete extinction of the disease in all forms and in all patches.  A partial extinction event is one in which one or more classes of infectious individuals goes extinct, but at least some class remains endemic.  Such events are transient from the prospective of deterministic and stochastic modeling and have not been considered to date.  Metapopulation models are characterized by multiple patches and the rates of movement between them.  It is of particular interest to consider partial extinction events in a metapopulation in which the disease goes extinct in some, but not all patches.  Statistics like the probability of partial extinction events may help to understand how the underlying structure of the metapopulation influences the dynamics of the system.  Additionally, the probability of extinction in a single patch of a metapopulation model may be viewed as a numerical rating of how susceptible that patch is to outbreak of disease.  When a patch corresponds to a locality, this rating could then be used to optimize control strategies from the perspective of that patch.    An attempt to study partial extinction events for an outbreak of ISAv in two patches using MTBP techniques led to the determination that these techniques are not suitable to answer such questions.

\section{two-patch model of ISAv}\label{2patch}
We begin by illustrating the use of MTBP to approximate the probability of extinction in metapopulation models by taking a two-patch model of ISAv as a test case.  The CTMC is constructed so that it is related to a previously studied deterministic model \cite{Milliken2016}.  As a result, the quasi-steady state is equal to the endemic equilibrium of the deterministic model.  Parameters are chosen to ensure the quasi-steady state associated to outbreak exists and can be easily located numerically.  They are also chosen to ensure the accuracy of the MTBP approximation.  They are not chosen for biological relevance.
\subsection{Deterministic SIV-SIV model}\label{det6}
In previous work with S. S. Pilyugin \cite{Milliken2016}, we proposed a two-patch SIV model to study the dynamics of an ISAv infection.  The two patches are coupled solely via diffusion of the virus.  Birth and death rates are patch dependent and are denoted by a subscript associated to the patch.  All other parameters are patch independent.  The force of infection in the $i^{th}$ patch is given by $S_i(\sigma I_i+\rho V_i)$, but the parameters $\sigma$ and $\rho$ can be scaled away.  Rescaling yields the following system:
\begin{equation}\label{6d}
\begin{aligned}\begin{cases}
			 \overset{.}{S_1}&=S_1(\beta_1-\mu_1 S_1)-S_1( I_1+ V_1)\\
			\overset{.}{I_1}&=S_1(I_1+ V_1)-\al I_1\\
			\overset{.}{V_1}&=k(V_2-V_1)-\om V_1+\delta I_1\\
			 \overset{.}{S_2}&=S_2(\beta_2-\mu_2 S_2)-S_2( I_2+ V_2)\\
			\overset{.}{I_2}&=S_2(I_2+ V_2)-\al I_2\\
			\overset{.}{V_2}&=k(V_1-V_2)-\om V_2+\delta I_2.
\end{cases}
\end{aligned}
\end{equation}
where $\beta_1,\,\beta_2$ are the patch specific birth rates of susceptible fish, $\mu_1,\,\mu_2$ are the patch specific, density dependent mortality rates, $\al$ is the mortality rate of infected fish, $\delta$ is the rate at which infected fish shed the virus into the environment, $\om$ is the rate at which it clears from the environment and $k$ is the rate of viral diffusion.

System \eqref{6d} admits 7 equilibria in total.  Four equilibria corresponding to the absence of the virus: (0,0,0,0,0,0), $(\bar{S}_1,0,0,0,0,0)$, $(0,0,0,\bar{S}_2,0,0)$, DFE $=(\bar{S}_1,0,0,\bar{S}_2,0,0)$.  Of these, only the disease-free equilibrium (DFE) is locally stable in the subspace associated to the absence of the disease.  Let 
$$\mR_0^{(1)}=\frac{(\om(2k+\om)+\delta(k+\om))\beta_1}{\al\om(2 k+\om)\mu_1}\quad\quad\quad\quad\mR_0^{(2)}=\frac{(\om(2k+\om)+\delta(k+\om))\beta_2}{\al\om(2 k+\om)\mu_2}.$$
Then $\mathcal{R}_0^{(i)}$ is the patch specific reproduction numbers corresponding to host fish only in patch $i$.  System \eqref{6d} admits two additional equilibria corresponding to the case where there are host fish only in patch one or only in patch two: $(S_1',I_1',V_1',0,0,V_2')\iff\mR_1^0>1$ and $(0,0,V_1^*,S_2^*,I_2^*V_2^*)\iff\mR_2^0>1$.  The basic reproduction number for system \eqref{6d} is given by
$$\mR_0=\hlf\left(\mR_1^0+\mR_2^0+\sqrt{(\mR_1^0-\mR_2^0)^2+4\bar{S}_1\bar{S}_2C^2}\right),$$
where $C=\frac{\delta k}{\al\om(2 k+\om)}.$  Following \cite{Milliken2016}, we have that DFE is globally asymptotically stable (\emph{g.a.s.}) if and only if $\mR_0\leq 1$.  If $\mR_0>1$, then the DFE is unstable and the virus invades and persists when introduced.  In fact, the subset of the boundary associated to the extinction of the virus is a uniform strong repeller whenever $\mR_0>1$ \cite{Butler1986,Fonda1988,Freedman1994,Garay1989,Hofbauer1989,Milliken2016,Thieme1993}.  If, in addition, the following symmetric conditions are met,
$$\mR_1^0>\frac{\mu_2}{\mu_1}Q(\mR_2^0-1)\quad\quad\text{and}\quad\quad\mR_2^0>\frac{\mu_1}{\mu_2}Q(\mR_1^0-1),$$
where $Q=\frac{\delta k}{\om(2k+\om)+\delta(k+\om)}$, then there exists a unique positive endemic equilibrium.

\subsection{Stochastic SIV-SIV model} 
From the preceding deterministic model we construct the CTMC, $\mathbf{X}(t)=(S_1(t),I_1(t),V_1(t),S_2(t),I_2(t),V_2(t))$, with the infinitesimal transition probability to state $j$ from state $i$ given by
$$p_{i,j}(\Delta t)=\mathbb{P}\{\mathbf{X}(t+\Delta t)=j\mid\mathbf{X}(t)=i\}=\sigma(i,j)\Delta t + o(\Delta t),$$
where $\sigma(i,j)$ is the rate associated to the transition from state $i$ to state $j$ and can be found in Table \ref{table:6drates}.
\begin{table}[H]
\begin{tabular}{l l c}\hline
Description&Transition& rate $\sigma(i,j)$\\\hline
Birth of $S_1$		&$(S_1,I_1,V_1,S_2,I_2,V_2)\mapsto(S_1+1,I_1,V_1,S_2,I_2,V_2)$	& $\beta_1 S_1$\\
Death of $S_1$		&$(S_1,I_1,V_1,S_2,I_2,V_2)\mapsto(S_1-1,I_1,V_1,S_2,I_2,V_2)$		&$\mu_1 S_1^2$\\
Infection of $S_1$	&$(S_1,I_1,V_1,S_2,I_2,V_2)\mapsto(S_1-1,I_1+1,V_1,S_2,I_2,V_2)$	&$S_1(I_1+V_1)$\\
Death of $I_1$		&$(S_1,I_1,V_1,S_2,I_2,V_2)\mapsto(S_1,I_1-1,V_1,S_2,I_2,V_2)$		&$\al I_1$\\
Shedding of $V_1$	&$(S_1,I_1,V_1,S_2,I_2,V_2)\mapsto(S_1,I_1,V_1+1,S_2,I_2,V_2)$	&$\delta I_1$\\
Clearance of $V_1$	&$(S_1,I_1,V_1,S_2,I_2,V_2)\mapsto(S_1,I_1,V_1-1,S_2,I_2,V_2)$		&$\om V_1$\\
Diffusion of $V_1$	&$(S_1,I_1,V_1,S_2,I_2,V_2)\mapsto(S_1,I_1,V_1-1,S_2,I_2,V_2+1)$	&$k V_1$\\
Birth of $S_2$		&$(S_1,I_1,V_1,S_2,I_2,V_2)\mapsto(S_1,I_1,V_1,S_2+1,I_2,V_2)$	& $\beta_2 S_2$\\
Death of $S_2$		&$(S_1,I_1,V_1,S_2,I_2,V_2)\mapsto(S_1,I_1,V_1,S_2-1,I_2,V_2)$		&$\mu_2 S_2^2$\\
Infection of $S_2$	&$(S_1,I_1,V_1,S_2,I_2,V_2)\mapsto(S_1,I_1,V_1,S_2-1,I_2+1,V_2)$	&$S_2(I_2+V_2)$\\
Death of $I_2$		&$(S_1,I_1,V_1,S_2,I_2,V_2)\mapsto(S_1,I_1,V_1,S_2,I_2-1,V_2)$		&$\al I_2$\\
Shedding of $V_2$	&$(S_1,I_1,V_1,S_2,I_2,V_2)\mapsto(S_1,I_1,V_1,S_2,I_2,V_2+1)$	&$\delta I_2$\\
Clearance of $V_2$	&$(S_1,I_1,V_1,S_2,I_2,V_2)\mapsto(S_1,I_1,V_1,S_2,I_2,V_2-1)$		&$\om V_2$\\
Diffusion of $V_2$	&$(S_1,I_1,V_1,S_2,I_2,V_2)\mapsto(S_1,I_1,V_1+1,S_2,I_2,V_2-1)$		&$k V_2$\\
\hline
\end{tabular}
\caption{State transitions and rates for the two-patch CTMC model, $X_t$.}
\label{table:6drates}
\end{table}
\begin{rmk} Recall that the original force of infection in the $i^{th}$ patch given by  $S_i(\sigma I_i+\rho V_i)$.  $S$ and $V$ are rescaled and $\mu$ and $\delta$ relabeled yielding \eqref{6d} for easier analysis.  The $V$ that is retained represents a scalar multiple of the number of virions present.  Let $\mu$ and $\delta$ reflect $\sigma=1$ and $\rho$ chosen so that we may interpret 1 unit of $V$ as any number of virions, such as an average infectious viral dose (e.g. ID50).  This makes the transition $V\mapsto V+1$ in the rescaled model reasonable.
\end{rmk}

We are interested in studying the dynamics after infectious agents are introduced to an entirely susceptible system.  Analysis of the flow of \eqref{6d} on the boundary shows that, in absence of the disease, DFE is \emph{g.a.s.}.  Therefore, we assume that DFE is the initial state of the system prior to introduction of the disease.  As $X_t$ evolves in time, $S_1(t)$ and $S_2(t)$ evolve along with all the other state components.   To formulate the MTBP we first pass to embedded discrete time Markov chain (DTMC), $X_n$.  Next, suppose that $S_1(n)\equiv\bar{S_1}$ and $S_2(n)\equiv\bar{S_2}$, the disease-free populations of susceptible fish in patches 1 and 2, respectively and that each individual gives birth independently of other individuals.  Let $Z_n=(I_1(n),V_1(n),I_2(n),V_2(n))$ be the random variable associated to the $n^{th}$ generation.  The offspring probability generating function (pgf) is given by 
$$\mathbf{F}(\mathbf{u})=(f_1(\mathbf{u}),f_2(\mathbf{u}),f_3(\mathbf{u}),f_4(\mathbf{u})),$$
where, for $i=1,2,3,4,$
$$f_i((u_1,u_2,u_3,u_4))=\sum_{n=0}^\infty p_i(r_1,\dotsc,r_4)u_1^{r_1}\dots u_4^{r_4},$$
and $p_i(r_1,\dotsc,r_4)$ is the probability that an object of type $i$ gives birth to $r_1$ offspring of type 1, $\dotsc,$ and $r_4$ offspring of type $4$.  
The offspring pgf for $I_1$ is 
$$f_1(\mathbf{u})=\frac{\al+\delta u_1u_2+\bar{S}_1u_1^2}{\al+\delta+\bar{S}_1},$$
the offspring pgf for $V_1$ is
$$f_2(\mathbf{u})=\frac{\om+k u_4+\bar{S}_1u_1u_2}{\om+k+\bar{S}_1},$$
the offspring pgf for $I_2$ is
$$f_3(\mathbf{u})=\frac{\al+\delta u_3u_4+\bar{S}_2u_3^2}{\al+\delta+\bar{S}_2},$$
and the offspring pgf for $V_2$ is 
$$f_4(\mathbf{u})=\frac{\om+k u_2+\bar{S}_2u_3u_4}{\om+k+\bar{S}_2}.$$
The matrix of expectations $\mathbb{M}=D\mathbf{F}(\mathbf{1})$ is given by
$$\mathbb{M}=\begin{bmatrix} \frac{\delta+2\bar{S}_1}{\al+\delta+\bar{S}_1}&\frac{\delta}{\al+\delta+\bar{S}_1}&0&0\\[5pt]
\frac{\bar{S}_1}{\om+k+\bar{S}_1}&\frac{\bar{S}_1}{\om+k+\bar{S}_1}&0&\frac{k}{\om+k+\bar{S}_1}\\[5pt]
0&0&\frac{\delta+2\bar{S}_2}{\al+\delta+\bar{S}_2}&\frac{\delta}{\al+\delta+\bar{S}_2}\\[5pt]
0&\frac{k}{\om+k+\bar{S}_2}&\frac{\bar{S}_2}{\om+k+\bar{S}_2}&\frac{\bar{S}_2}{\om+k+\bar{S}_2}\\[5pt]\end{bmatrix}.$$

A branching process is called positively regular if $\mathbb{M}$ is primitive.  A $k$--many type process is called not singular if $F(\mathbf{0})>\mathbf{0}$ with respect to the standard order, and whenever $\mathbf{x},\mathbf{y}\in[0,1]^k$ with $\mathbf{x}\leq\mathbf{y}$, then $D\mathbf{F}(\mathbf{x})\leq D\mathbf{F}(\mathbf{y})$.  The $i,j^{th}$ entry of $D\mathbf{F}(\mathbf{1})$ is $\frac{\partial f_i}{\partial u_j}(\mathbf{1})$, the expected number of type $j$ offspring of an individual of type $i$.  Following Harris \cite{Harris1963}, let $q_i$ be the extinction probability if initially there is one object of type $i$, $i=1,\dotsc,k$.  Let $\mathbf{q}=(q_1,\dotsc,q_k)$.  Let $\mathbb{P}_0$ be the probability of extinction of the branching process given that $Z_0=(j_1,\cdots,j_k)$.  Since we have assumed that individuals give birth independent of one another, 
$$\mathbb{P}_0=q_1^{j_1}q_2^{j_2}\dots q_k^{j_k}.$$
  
The branching process constructed above to approximate ISAv in two patches is positively regular (in fact, $\mathbb{M}^3>0$).  It is easily verified that it is also not singular.  
The Threshold Theorem of Allen and van den Driessche \cite{Allen2013} and Theorem 7.1 of Harris \cite{Harris1963} combine to give the following result.
\begin{thm}\label{crit} Suppose $Z_n$ is a MTBP with probability generating function $\mathbf{F}:\mathbb{R}^k\to\mathbb{R}^k$ such that $\mathbf{F}(\mathbf{0})>\mathbf{0}$, $D\mathbf{F}(\mathbf{x})\leq D\mathbf{F}(\mathbf{y})$ whenever $\mathbf{x},\mathbf{y}\in[0,1]^k$ with $\mathbf{x}\leq\mathbf{y}$, and $D\mathbf{F}(\mathbf{1})$ is primitive.  If $\mathcal{R}_0\leq1$, then $\mathbf{q}=\mathbf{1}$.  If $\mathcal{R}_0>1$, then $\mathbf{q}$ is the unique vector $\mathbf{0}\leq\mathbf{q}<\mathbf{1}$ satisfying $\mathbf{F}(\mathbf{q})=\mathbf{q}$. 
\end{thm}
\begin{rmk} For a fixed initial vector $z_{_0}$, the probability of extinction $\mathbb{P}_0=P(Z_n=\mathbf{0}|Z_0=z_{_0}\text{ for some }n>0)$.  In a metapopulation model, the probability of extinction is, therefore, the probability that all infectious classes go extinct, in all patches.  
If we wanted to use MTBP approximation to calculate a partial extinction event, like extinction in one patch, we would have to recast the MTBP to only track the evolution of those infectious classes and assume the number of individuals in other infectious classes remain fixed.  However, we already assumed that there are few individuals initially present in each infectious class.  As we have discussed above, in order to justify the assumption that the number of individuals in a given class remains fixed, the initial population in that class must be sufficiently large.  The MTBP is, therefore, not the appropriate tool to study partial extinction events.
\end{rmk}

\subsection{Numerical example}
In order to illustrate the accuracy of MTBP approximation of the probability of total extinction in a metapopulation model, we choose parameter values according to two criteria: $(i)$ the disease-free number of susceptible fish is sufficiently large in each patch for approximation by branching process; and $(ii)$ the endemic equilibrium of the deterministic system \eqref{6d} can be located numerically.  The endemic equilibrium of \eqref{6d} is a quasi-steady state of the CTMC and the embedded DTMC.  The second criterion also implies that $\mathcal{R}_0>1$. 
 Therefore, purely for the purpose of illustration and without regard to biological relevance, we consider the parameter vector $(\beta_1=4,\mu_1=0.05,\beta_2=2.4,\mu_2=0.04,\al=3.3,\delta=1.3,\om=4,k=3)$.  Then $\bar{S}_1=80$, $\bar{S}_2=60$, $\mR_1^0\approx 30$, $\mR_2^0\approx 22$, and $\mR_0\approx30>>1$.  Recall that $\mathbb{P}_0=q_1^{j_1}q_2^{j_2}q_3^{j_3}q_4^{j_4}$, where $Z_0=(j_1,j_2,j_3,j_4)$ and $q_i$ is the extinction probability if there is initial one object of type $i$.  Because of this and due to the computational expense of simulating this model, we only consider initial states with one object of type $i$, $i=1,\dotsc,4$.  The vector $\mathbf{q}$ of extinction probabilities is determined by iterating the pgf from the initial vector $\mathbf{0}$.  Let $\mathbb{P}_0^{(n)}$ denote the probability of extinction approximated by numerical simulation over $n$ realizations.  The results are presented in Table \ref{table:num6d}.
\begin{table}[h!]
\begin{tabular}{l l l l c r}\hline
$I_1(0)$&$V_1(0)$&$I_2(0)$&$V_2(0)$&$\mathbb{P}_0$& $\mathbb{P}_0^{(1,000,000)}$\\\hline
1	&0&0&0	&	0.0406	&0.0410\\
0	&1&0&0	&	0.0501	&0.0501\\
0	&0&1&0	&	0.0538	&0.0542\\
0	&0&0&1	&	0.0650	&0.0652\\
\hline
\end{tabular}
\caption{Probability of extinction of the virus from the initial state $(\bar{S}_1,I_1(0),V_1(0),\bar{S}_2,I_2(0),V_2(0))$ and parameter vector $(\beta_1=4,\mu_1=0.05,\beta_2=10,\mu_2=0.04,\al=3.3,\delta=1.3,\om=4, k=3)$ is approximated by MTBP and numerically over $1,000,000$ realizations.}
\label{table:num6d}
\end{table}

By the law of large numbers, as the number of realizations, $n$, increases to infinity, $\mathbb{P}_0^{(n)}$ tends to the true probability of extinction.  Assuming that $\mathbb{P}_0^{(n)}$ is distributed normally, the error in approximating $\mathbb{P}_0$ with $\mathbb{P}_0^{(n)}$ goes to zero like $\frac{1}{\sqrt{n}}$.  This implies that approximation of $\mathbb{P}_0$ to three decimal places by numerical simulation requires making $10^6$ realizations, at great computational expense.

The results in Table \ref{table:num6d} suggest that the MTBP approximates the probability of extinction in the CTMC very accurately.  In this case, we were not able to solve the nonlinear system of equations given by $\mathbf{F}(\mathbf{q})=\mathbf{q}$ for an analytical solution to the MTBP.  However, we are able to approximate $\mathbf{q}$ by iteration with little computational expense, since $\mathbf{F}^n(\mathbf{0})\to\mathbf{q}$.

\section{one-patch model}\label{1patch}
As discussed above, we must assume the number of susceptible individuals remains fixed at the disease-free level in order to utilize branching process techniques for SIV models.  The disease-free population size must be at least as large as some critical value in order for this assumption to be reasonable.  Currently, there is no analytic estimate of this critical size.  In this section and section \ref{size}, we compare MTBP approximation and simulation of the CTMC at a range of small initial populations for two models.  These models are introduced in this section and the next.  The first is an invariant subsystem of \eqref{6d}, which models the dynamics of infection in a single patch.  We consider this one-patch model because it reduces the computational expense while still retaining the key features of interest.
\subsection{Deterministic SIV model}

When there is no diffusion, i.e. $k=0$, then each patch of the two-patch system forms an invariant SIV subsystem given by:
\begin{equation}\label{coup1}
\begin{aligned}\begin{cases}
			 \overset{.}{S}&=S(\beta-\mu S)-Sf(I,V)\\
			\overset{.}{I}&=Sf(I,V)-\al I\\
			\overset{.}{V}&=-\om V+\delta I.
\end{cases}
\end{aligned}
\end{equation}
where $f(I,V)=f_1(I,V)=(I+V)$, $\beta$ is the birth rate of susceptible fish, $\mu$ the mortality rate of susceptible fish, $\al$ the mortality rate of infected fish, $\om$ is the rate of viral clearing and $\delta$ is the rate of viral shedding. All of these parameters are assumed to be positive.\p

The system admits equilibria (0,0,0) (which is always unstable) and the disease-free equilibrium (DFE), ($\bar{S}$,0,0).  The basic reproduction number is, 
\begin{equation}\label{MAR0}\mR_0=\frac{(\delta+\om)\beta}{\al\om\mu}.\end{equation}
When $\mR_0>1$ the system also admits a unique positive endemic equilibrium.  $\mR_0=1$ is also a threshold for the dynamics of the system.  If $\mR_0\leq1$, then the DFE is \emph{g.a.s.}.  If $\mR_0>1$, then the DFE is unstable and the virus invades and persists when introduced.  The largest invariant subset of the boundary is a uniform strong repeller when $\mR_0>1$ \cite{Milliken2016,Thieme1993}.

\subsection{Stochastic SIV model}\label{sec:stoch3d}

The CTMC model $\mathbf{X}(t)=(S(t),I(t),V(t))$ associated to system \eqref{coup1} with $f(I,V)=f_1(I,V)$ is characterized by the transition rates given in Table \ref{table:3drates}.
\begin{table}[h!]
\begin{tabular}{l l c}\hline
Description&Transition& rate $\sigma(i,j)$\\\hline
Birth of $S$&$(S,I,V)\mapsto(S+1,I,I)$& $\beta S$\\
Death of $S$& $(S,I,V)\mapsto(S-1,I,V)$&$\mu S^2$\\
Infection&$(S,I,V)\mapsto(S-1,I+1,V)$&$S(I+V)$\\
Death of $I$&$(S,I,V)\mapsto(S,I-1,V)$&$\al I$\\
Shedding of $V$&$(S,I,V)\mapsto(S,I,V+1)$&$\delta I$\\
Clearance of $V$&$(S,I,V)\mapsto(S,I,V-1)$&$\om V$\\\hline
\end{tabular}
\caption{State transitions and rates for the CTMC SIV model.}
\label{table:3drates}
\end{table}
\begin{figure}
\includegraphics[width=\textwidth]{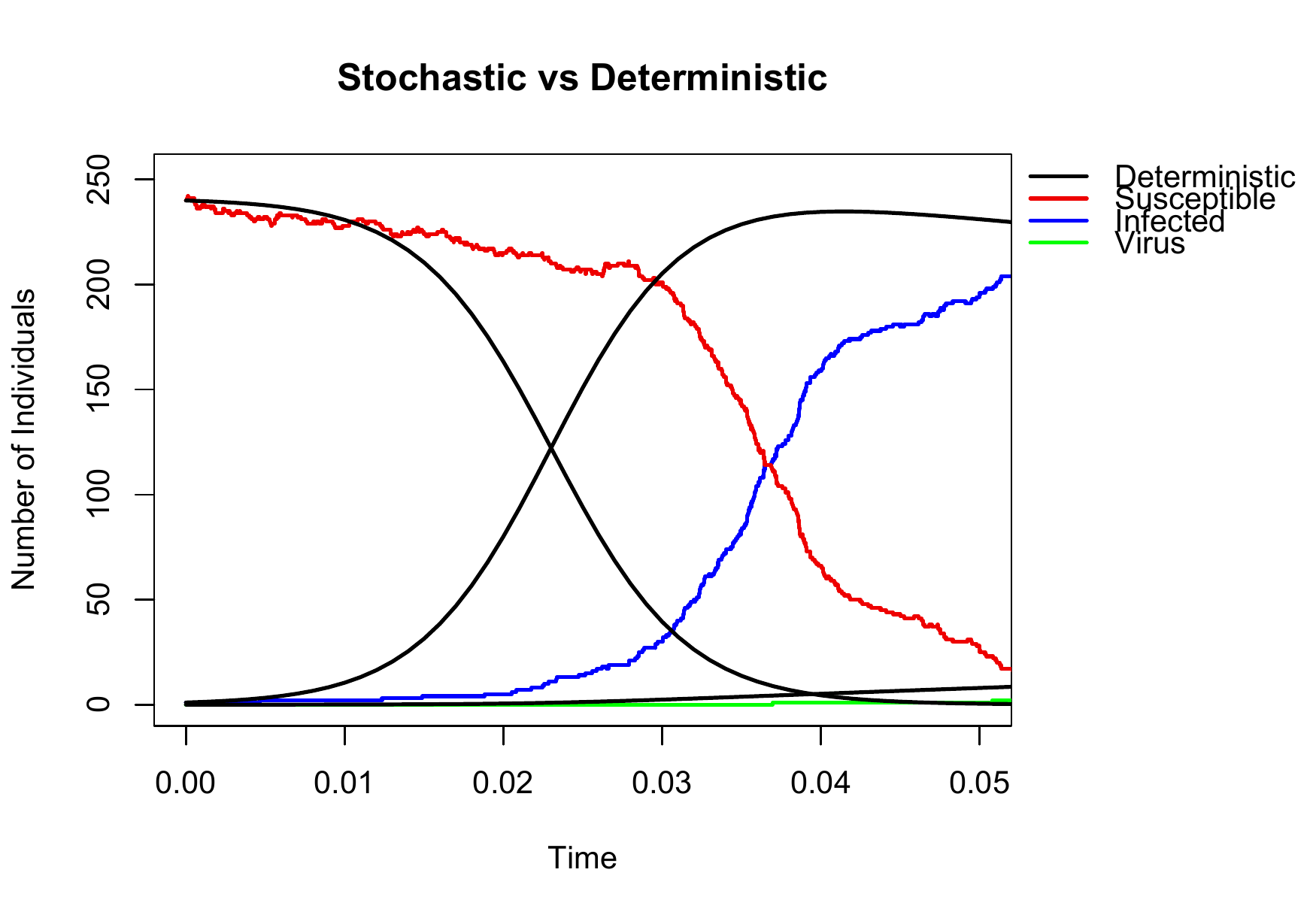}
\caption{One realization of the CTMC model, $X_t$, compared to solution of the deterministic model.  Both simulations take initial condition $(S=240,I=1,V=0)$ and parameter vector $(\beta=12,\mu=0.05,\alpha=3.3,\delta=1.3,\om=4)$.}
\end{figure}

To estimate the probability of extinction of the virus, we approximate the CTMC near the DFE.  As in the two-patch case, we pass to the embedded DTMC, assume that $S(n)\equiv \bar{S}$ and that all individuals give birth independently.  Let $Z_n=(I(n),V(n))$ and construct the probability-generating function (pgf) for the MTBP, $Z_n$.  
$$\mathbf{F}(\mathbf{u})=\left(f_1(\mathbf{u}),f_2(\mathbf{u})\right)=\left(\frac{\al+\delta u_1u_2+\bar{S}u_1^2}{\al+\delta+\bar{S}},\frac{\om+\bar{S}u_1u_2}{\om+\bar{S}}\right).$$
It follows that $Z_n$ is not singular and the matrix of expectations is given by
$$\mathbb{M}=\begin{bmatrix} \frac{\delta+2\bar{S}}{\al+\delta+\bar{S}}&\frac{\delta}{\al+\delta+\bar{S}}\\[5pt]
\frac{\bar{S}}{\om+\bar{S}}&\frac{\bar{S}}{\om+\bar{S}}\\[5pt]
\end{bmatrix},$$
is positive.  Thus, $\mathbb{M}$ is primitive and Theorem \ref{crit} applies.  Solving the system of nonlinear equations given by $\mathbf{F}(\mathbf{q})=\mathbf{q}$ yields
\begin{equation}\label{MAq1}
q_1=\frac{\al+\delta+\om+\bar{S}-\sqrt{(\al-(\om+\bar{S}))^2+\delta(\delta+2(\al+\om+\bar{S})}}{2\bar{S}},\text{ and}
\end{equation}
\begin{equation}q_2=\frac{\om}{\om+\bar{S}(1-q_1)}.\end{equation}
Then the probability of extinction given that $Z_0=(j_1,j_2)$ is
 $$\mathbb{P}_0=q_1^{j_1}q_2^{j_2}.$$ 
Note that, for this model, the MTBP approximation of the probability of extinction can be determined analytically.  That is, $\mathbb{P}_0$ can be expressed as a continuous function of the parameters.
\subsection{Numerical example}
For the purpose of illustrating the accuracy of the MTBP approximation, we consider the parameter vector given by $(\beta=4,\mu=0.05,\al=3.3,\delta=1.3,\om=4)$.  This choice of parameters yields $\bar{S}=80$ and $\mR_0\approx32>>1$.  Let $\mathbb{P}_0$ denote the probability of extinction predicted by the MTBP, given $Z_0=(I(0),V(0))$.  The probability of extinction in the CTMC is estimated by simulating numerically.  Let $\mathbb{P}_0^{(1,000,000)}$ denote the probability of extinction approximated by numerical simulation over $1,000,000$ realizations.  The results of both approximations are presented in Table \ref{table:num3d}.

\begin{table}[H]
\begin{tabular}{l l c c c r}\hline
$I(0)$&$V(0)$&$\mathbb{P}_0$& $\mathbb{P}_0^{(1,000,000)}$\\\hline
1	&0	&	0.0406		&0.0407\\
0	&1	&	0.0495		&0.0494\\
1	&1	&	0.0020		&0.0020\\\hline
\end{tabular}
\caption{Probability of extinction of the virus from the initial condition $(\bar{S},I(0),V(0))$ with the parameter vector $(\beta=4,\mu=0.05,\al=3.3,\delta=1.3,\om=4)$ approximated by branching process and numerically over $1,000,000$ realizations.}
\label{table:num3d}
\end{table}
\begin{figure}[H]
  \includegraphics[width=\linewidth]{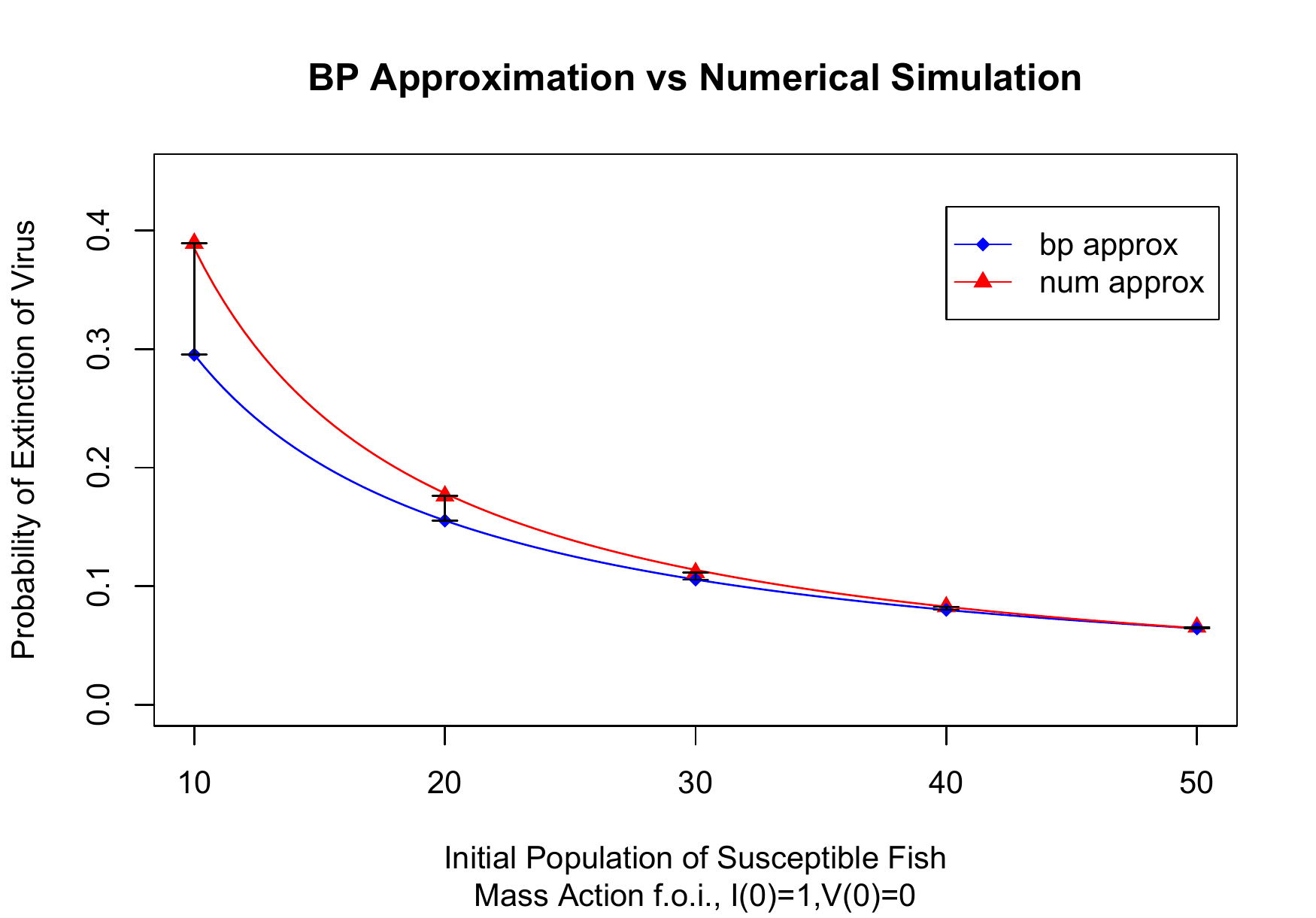}
  \caption{Comparison of multitype branching process approximation to numerical simulation of probability of extinction in single patch model with mass action force of infection and high mortality for infected fish.}\label{fig:MAcrit}
  \label{MABeta}
\end{figure}

Table \ref{table:num3d} illustrates that, for this choice of parameters, the MTBP provides extremely accurate results.  Since $\mathbb{P}_0$ can be expressed as a function of the parameters, the computational expense for MTBP approximation is negligible.  However, we cannot be certain, a priori, whether or not the disease-free population of susceptible fish is sufficiently large without comparing the MTBP results to numerical simulation.  Therefore, the estimate of computational expense for MTBP approximation should include the cost of simulating the CTMC.   The additional expense for simulating the CTMC can be significant.  

In Figure \ref{MABeta}, we compare MTBP and numerical simulation for initial populations at ten unit increments from 10 to 50.  First, note that the population of susceptible fish at DFE is given by $\bar{S}=\frac{\beta}{\mu}$.  Therefore, by assuming $\mu=1$, we have that $\bar{S}=\beta$.  We fix the remaining parameters $(\mu=1,\al=3.3,\delta=1.3,\om=4)$ and vary $\beta$ from 10 to 50 in ten unit increments.

  Numerical data in Figure \ref{fig:MAcrit} is fit with a power law curve $y=bx^\lambda$ where $b=4.9584$ and $\lambda=-1.11$.  Not pictured, the absolute error is fit with a power law curve with $b=62.172$ and $\lambda=-2.743$ and the relative error is fit with a power law curve with $b=0.4584$ and $\lambda=-0.067$.  Since $\mathbb{P}_0$ is a continuous function of the parameters, there was no need to fit a curve to the MTBP results.
  
In Section \ref{size}, we will show that the character and speed of convergence of the MTBP approximation results to the CTMC simulation results depends on the structure of the model and the choice of parameters.  We do this by constructing illustrations similar to Figure \ref{MABeta} based on variations of the one-patch model.  

\section{one-patch model with modified force of infection}\label{MM}
\subsection{Deterministic model}
The one-patch model given by system \eqref{coup1} proposes a mass action force of infection.  It has been suggested that the \emph{f.o.i.} may initially be driven by infected salmon encountering susceptible salmon when there are low levels of free virus present at the outset of an exposure event.  As more salmon become infected and shed more and more virus into the environment, the free virus may then drive the infection.  To account for this we modify system \eqref{coup1} by considering $f(I,V)=f_2(I,V)$ where 
$$f_2(I,V)=\frac{m_1I}{a_1+I+V}+\frac{m_2V}{a_2+I+V}$$
Note that when $m_1=m_2$ and $a_1=a_2$, the growth function $S\left(\frac{m_1I}{a_1+I+V}+\frac{m_2V}{a_2+I+V}\right)$ simplifies to the standard Michelis-Menten function for $I+V$.  System \eqref{coup1} with $f(I,V)=f_2(I,V)$ admits equilibria at $\mathbf{0}$ and the DFE $(\frac{\beta}{\mu},0,0)$.  Following the next generation matrix approach \cite{Diekmann1990,VanDenDriessche2002}  the basic reproduction number is determined to be
\begin{equation}\label{MonR0}\mR_0=\frac{m_1a_2+\frac{\delta}{\om}m_2a_1}{\al a_1a_2}\frac{\beta}{\mu}.\end{equation}
The endemic equilibrium is a root of the vector field.  From $\overset{.}{V}=0$ we have $V'=\frac{\delta}{\om}I'$.  Substituting into $\overset{.}{I}=0$ yields $S'=f_1(I')$.  Let $f_2(I')=\frac{m_1I}{a_1+(1+\frac{\delta}{\om})I'}$ and  $f_3(I')=\frac{m_2\frac{\delta}{\om}I'}{a_2+(1+\frac{\delta}{\om})I'}$.  Then the nonnegative root of $\overset{.}{S}=0$ is a root of the equation
\begin{equation}\label{eequil}\beta-\mu\al f_1(I')-f_2(I')-f_3(I')=0.\end{equation}
Furthermore, $f_1'(I'),\;f_2'(I'),\;f_3'(I')>0$ and $f_2(0)=f_3(0)=0$.  Thus, \eqref{eequil} has a unique positive root if and only if $f_1(0)<\beta\iff\mR_0>1$.  Thus, the unique positive endemic equilibrium exists if and only if $\mR_0>1$.  If $\mR_0\leq 1$, then the DFE is \emph{g.a.s.}.  This system has the same dynamics on the boundary as the system with mass action \emph{f.o.i.}.  Using arguments similar to those in \cite{Milliken2016}, it follows that system \eqref{coup1} with $f(I,V)=f_2(I,V)$ is uniformly strongly persistent whenever $\mathcal{R}_0>1$ \cite{Thieme1993}.
\subsection{Stochastic model}
The CTMC model related to system \eqref{coup1} with $f(I,V)=f_2(I,V)$ is characterized by the transitions and rates given in Table \ref{table:3dMrates}.
\begin{table}[h!]
\begin{tabular}{l l c}\hline
Description&Transition& rate $\sigma(i,j)$\\\hline
Birth of $S$&$(S,I,V)\mapsto(S+1,I,I)$& $\beta S$\\
Death of $S$& $(S,I,V)\mapsto(S-1,I,V)$&$\mu S^2$\\
Infection&$(S,I,V)\mapsto(S-1,I+1,V)$&$S(\frac{m_1I}{a_1+I+V}+\frac{m_2V}{a_2+I+V})$\\
Death of $I$&$(S,I,V)\mapsto(S,I-1,V)$&$\al I$\\
Shedding of $V$&$(S,I,V)\mapsto(S,I,V+1)$&$\delta I$\\
Clearance of $V$&$(S,I,V)\mapsto(S,I,V-1)$&$\om V$\\\hline
\end{tabular}
\caption{State transitions and rates for the CTMC SIV model.}
\label{table:3dMrates}
\end{table}
\begin{figure}[H]
\includegraphics[width=\textwidth]{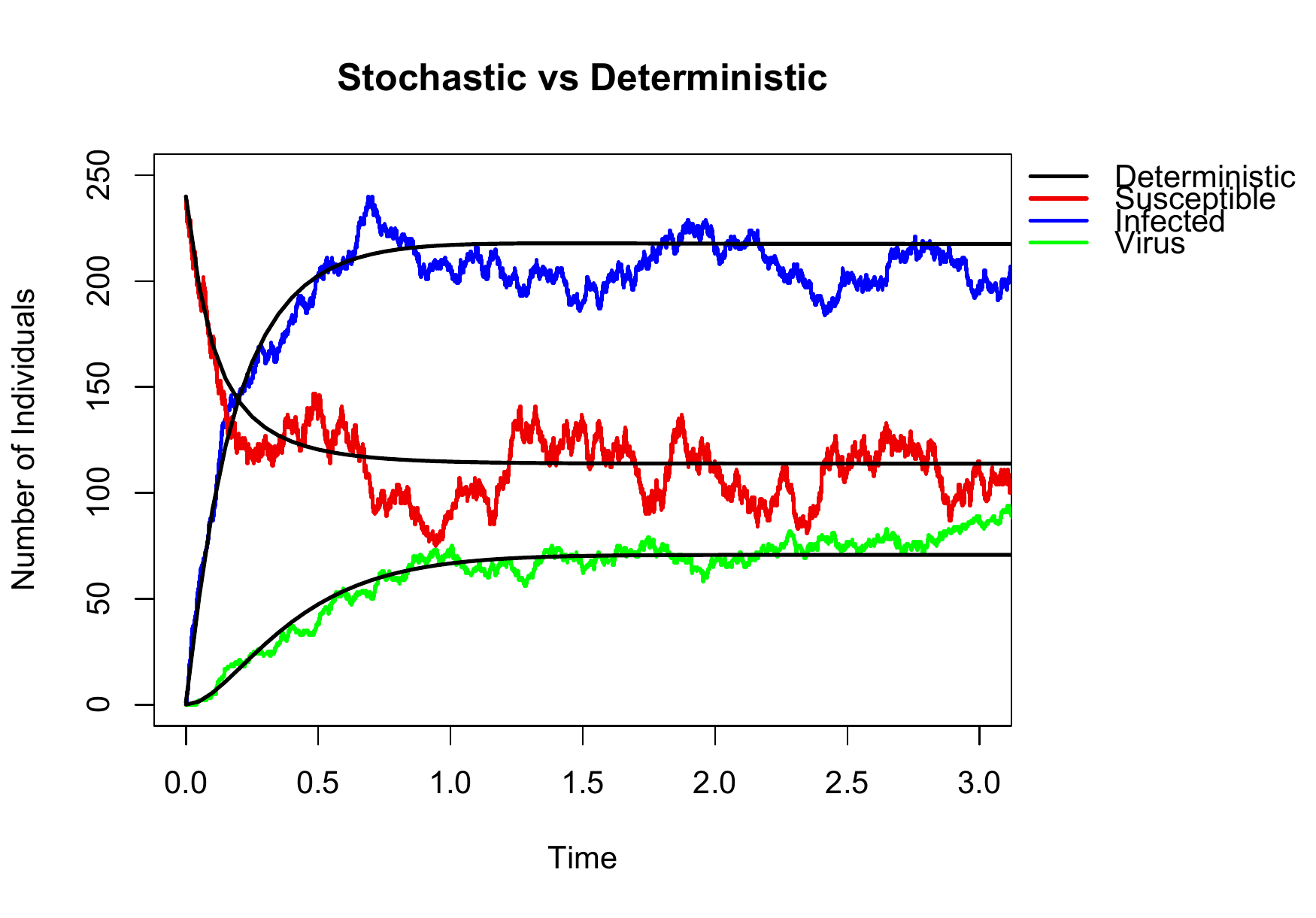}
\caption{One realization of the Markov chain model compared to solution of the deterministic model.  Both simulations take initial condition $(S=240,I=1,V=0)$ and parameter vector $(\beta=12,\mu=0.05,\alpha=3.3,\om=4,\delta=1.3,m_1=6,m_2=7.5,a_1=3,a_2=2)$.}
\end{figure}

We approximate the CTMC, $X_n$, near the DFE with the MTBP, $Z_n$, with the pgf 
$$\mathbf{F}(\mathbf{u})=\left(f_1(\mathbf{u}),f_2(\mathbf{u})\right)
=\left(
\frac{\al+\delta u_1u_2+ \bar{S}\frac{m_1}{a_1+1}u_1^2}{\al+\delta+\bar{S}\frac{m_1}{a_1+1}},
\frac{\om+\bar{S}\frac{m_2}{a_2+1}u_1u_2}{\om+\bar{S}\frac{m_2}{a_2+1}}\right).$$
The matrix of expectations is given by 
$$\mathbb{M}=\begin{bmatrix} \frac{\delta+2\bar{S}\frac{m_1}{a_1+1}}{\al+\delta+\bar{S}\frac{m_1}{a_1+1}}&\frac{\delta}{\al+\delta+\bar{S}\frac{m_1}{a_1+1}}\\[5pt]
\frac{\bar{S}\frac{m_2}{a_2+1}}{\om+\bar{S}\frac{m_2}{a_2+1}}&\frac{\bar{S}\frac{m_2}{a_2+1}}{\om+\bar{S}\frac{m_2}{a_2+1}}\\[5pt]\end{bmatrix}.$$
Clearly, the branching process in not singular and $\mathbb{M}$ is a positive matrix.  Thus, Theorem \ref{crit} applies.  Let $\Delta_1=\frac{m_1}{a_1+1}$, $\Delta_2=\frac{m_2}{a_2+1}$, and
\begin{equation}\mathscr{D}=\big(\al\Delta_2-\Delta_1(\bar{S}\Delta_2+\om)\big)^2+\delta\Delta_2^2(\delta+2\al+2\bar{S}\Delta_1)+2\delta\om\Delta_1\Delta_2.\end{equation}
  Then $\mathscr{D}>0$,
\begin{equation}\label{Monq1}q_1=\frac{\al\Delta_2+\delta\Delta_2+\om\Delta_1+\bar{S}\Delta_1\Delta_2-\sqrt{\mathscr{D}}}{2\bar{S}\Delta_1\Delta_2},\text{ and}\end{equation}
\begin{equation}\label{Monq2}q_2=\frac{\om}{\om+\bar{S}\Delta_2(1-q_1)}.\end{equation}
Given $Z_0=(j_1,j_2)$, $\mathbb{P}_0=q_1^{j_1}q_2^{j_2}$ can be expressed as a continuous function of the parameters.

\subsection{Numerical example}\label{parm}
For the purpose of illustrating the accuracy of the MTBP approximation, we consider the parameter vector given by $(\beta=4, \mu=0.05,\al=3.3,\om=4,\delta=1.3,m_1=3,m_2=2.5,a_1=3,a_2=2)$.  This implies that $\bar{S}=80$ and $\mR_0\approx34>>1$.  Let $\mathbb{P}_0$ denote the probability of extinction predicted by the MTBP, given $Z_0=(I(0),V(0))$.  The probability of extinction in the CTMC is estimated by simulating numerically.  Let $\mathbb{P}_0^{(1,000,000)}$ denote the probability of extinction approximated by numerical simulation over $1,000,000$ realizations.  The extinction probability predicted by the branching process approximation is compared with numerical results in Table \ref{table:num:3dM}.
\begin{table}[H]
\begin{tabular}{l l c c c r}\hline
$I(0)$&$V(0)$&$\mathbb{P}_0$&$\mathbb{P}_0^{(1,000,000)}$\\\hline
1	&0	&	0.0538		&0.0548\\
0	&1	&	0.0596		&0.0606\\
1	&1	&	0.0032		&0.0042\\\hline
\end{tabular}
\caption{Probability of extinction of the virus from the initial condition $(\bar{S},i_0,v_0)$ with the parameter vector $(\beta=4,\mu=0.05,\al=3.3,\om=4,\delta=1.3,m_1=3,m_2=2.5,a_1=3,a_2=2)$ approximated by branching process and numerically over $1,000,000$ realizations.}
\label{table:num:3dM}
\end{table}

This model represents a variant to the one-patch model studied in Section \ref{1patch} which differs only in the choice of function for the force of infection.
\section{Critical Size of disease-free Population}\label{size}
\begin{figure}[H]
  \includegraphics[width=\linewidth]{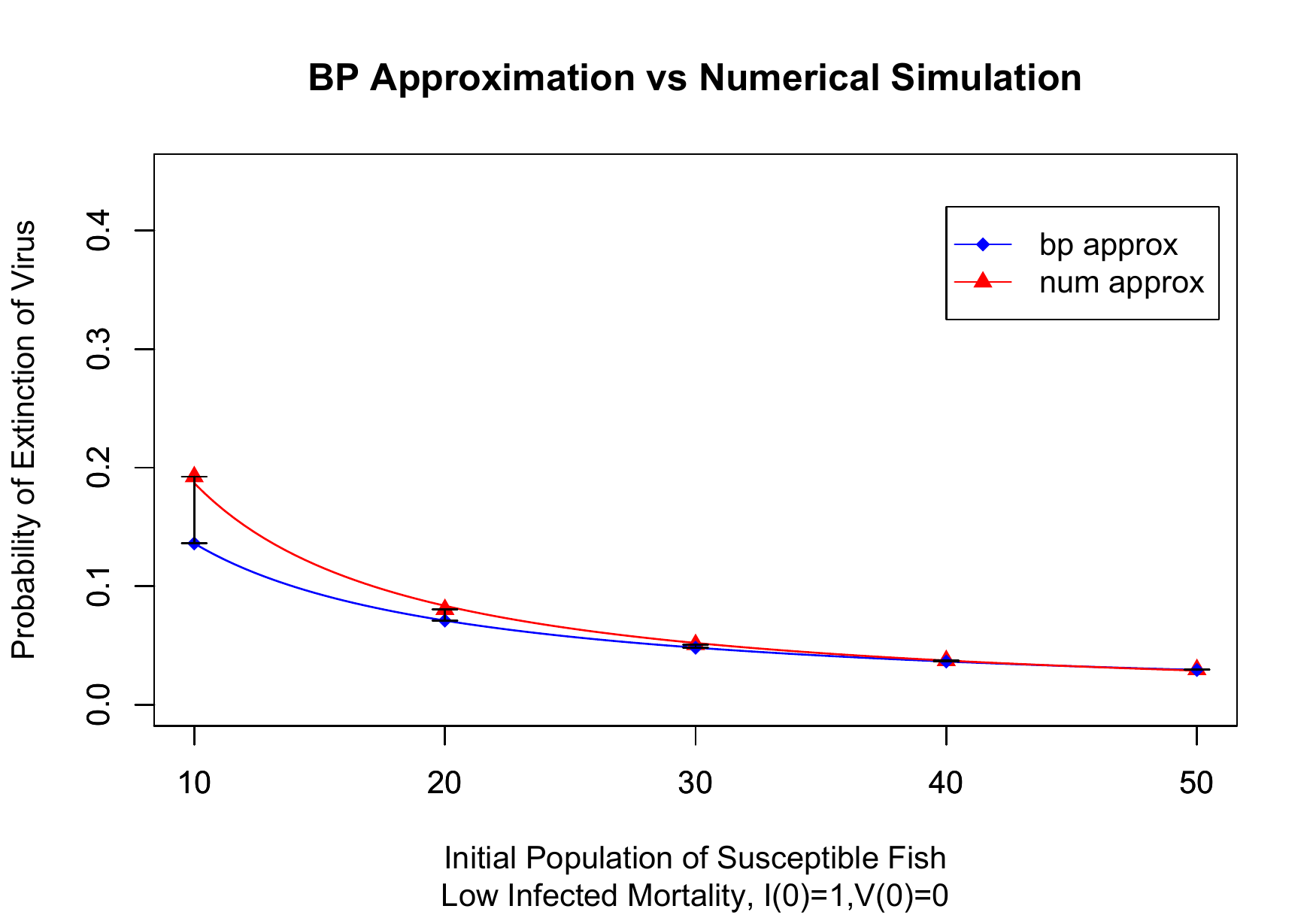}
    \caption{Comparison of multitype branching process approximation to numerical simulation of probability of extinction in single patch model with $f(I,V)=f_1(I,V)$ and parameter vector $(\mu=1,\al=1.5,\delta=1.3,\om=4)$.  }\label{IdBeta}
\end{figure}

In this section, we illustrate how variations to the underlying model affect the accuracy of MTBP approximation for small initial populations.  Figure \ref{MABeta} at the end of Section \ref{1patch} shows how MTBP approximation diverges from the probability of extinction in the CTMC for small initial populations when $f(I,V)=f_1(I,V)$.  We take this illustration as a baseline and vary the system in two ways.  First, we leave $f(I,V)=f_1(I,V)$, but reduce the mortality rate of infected fish from $\al=3.3$ to $\al=1.5$.  Second, we let $\al=3.3$ as in the baseline, but let $f(I,V)=f_2(I,V)$ as in the model developed in Section \ref{MM}.  In Figure \ref{IdBeta}, we set $f(I,V)=f_1(I,V)$ and fix the parameter vector $(\mu=1,\al=1.5,\delta=1.3,\om=4)$ with low mortality of infected fish and vary $\beta$ from 10 to 50 in ten unit increments.  Numerical data is fit with a power law curve $y=bx^\lambda$ where $b=2.7264$ and $\lambda=-1.164$.  Not pictured, the absolute error is fit with a power law curve with $b=42.109$ and $\lambda=-2.847$ and the relative error is fit with a power law curve with $b=14.459$ and $\lambda=-1.659$.
\begin{figure}
  \includegraphics[width=\linewidth]{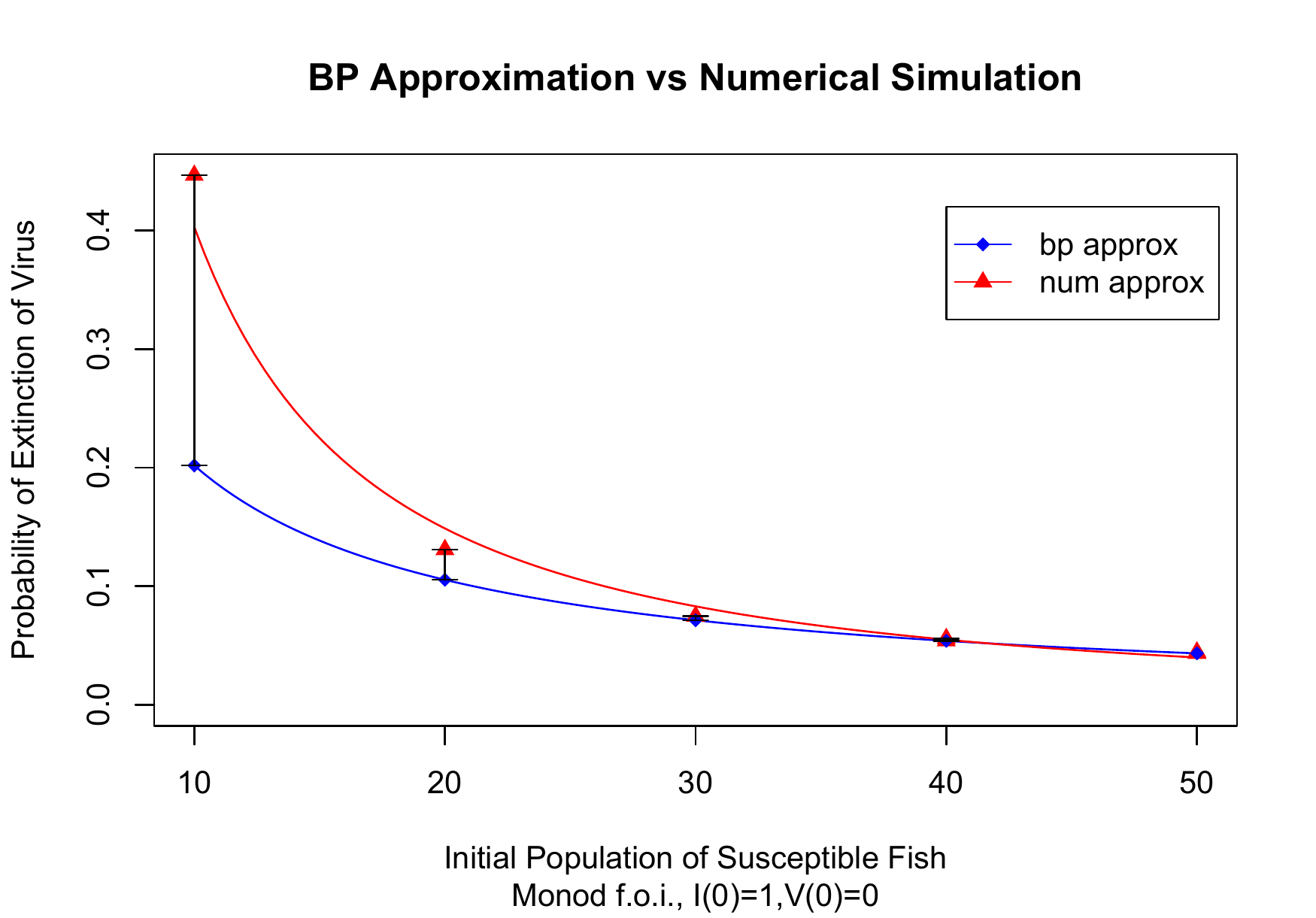}
     \caption{Comparison of multitype branching process approximation to numerical simulation of probability of extinction in single patch model with $f(I,V)=f_2(I,V)$ and parameter vector $(\mu=1,\al=3.3,\delta=1.3,\om=4,m_1=6,m_2=7.5,a_1=3,a_2=2)$.}\label{MonBeta}
     \end{figure}

In Figure \ref{MonBeta}, we set $f(I,V)=f_2(I,V)$ and fix the parameter vector $(\mu=1,\al=3.3,\delta=1.3,\om=4,m_1=6,m_2=7.5,a_1=3,a_2=2)$ and vary $\beta$ from 10 to 50 in ten unit increments.  Numerical data is fit with a power law curve $y=bx^\lambda$ where $b=11.074$ and $\lambda=-1.439$.  Not pictured, the absolute error is fit with a power law curve with $b=756.67$ and $\lambda=-3.507$ and the relative error is fit with a power law curve with $b=68.329$ and $\lambda=-2.068$.

Note that the results in Figures \ref{MABeta}, \ref{IdBeta}, and \ref{MonBeta} are graphed on the same axes on the same scale.  It is immediately evident that the character of convergence of the MTBP varies from the baseline illustration in each of the two latter ones.  It is harder to see from the graphs themselves, but the speed of convergence varies slightly as well.  This can be seen in Table \ref{Table:comp}.
\begin{table}[H]
\begin{tabular}{l r r r}\hline
Init. Pop.& $f_1,\alpha=3.3$&$f_1,\alpha=1.5$&$f_2,\alpha=3.3$\\\hline
10	&0.094	&0.056	&0.245			\\
20	&0.021	&0.009	&0.025			\\
30	&0.006	&0.003	&0.003			\\
40	&0.003	&0.001	&0.002			\\
50	&0.001	&0.000	&0.001			\\\hline
\end{tabular}
\caption{Entries represent absolute error between numerical results and multitype branching process results, $|\mathbb{P}_0-\mathbb{P}_0^{(1,000,000)}|$.  Column 2 corresponds to Figure \ref{MABeta}, Column 3 to Figure \ref{IdBeta} and Column 4 to Figure \ref{MonBeta}.}
\label{Table:comp}
\end{table}

\section{Discussion}  In this article, we use a model of ISAv in two patches and an invariant subsystem corresponding to one patch as toy models to develop CTMC models and MTBP approximations to estimate the probability of disease outbreak.  In addition to these models, we formulate a new one-patch model by varying the force of infection function.  In the case of the two-patch model, we approximate the probability of disease extinction, $\mathbb{P}_0$, by iterating the probability generating function of the MTBP.  For each one-patch model, characterized by its force of infection, it is possible to write the MTBP approximation of $\mathbb{P}_0$ as a continuous function of the parameters.  By comparing MTBP results to numerical simulation of the related CTMC, we show that, for large initial populations of susceptible fish, the MTBP approximation provides a good estimate of $\mathbb{P}_0$.  However, we should also note that MTBP approximation fails to provide accurate estimates of $\mathbb{P}_0$ when the initial population of susceptible fish is low.  It is therefore necessary to approximate $\mathbb{P}_0$ by numerical simulation concurrent with MTBP approximation.  While the computational expense for MTBP approximation is negligible, the computational expense for numerical simulation of the related CTMC, can be very high, especially for metapopulation models.  

In this article, we have not provided an analytical estimate on how large the initial population of susceptible individuals needs to be in order for the MTBP approximation to provide a good estimate of $\mathbb{P}_0$.  We have, however, illustrated the manner in which the approximation diverges from the true probability in several test cases.  Comparison of results in Figures \ref{MABeta}, \ref{IdBeta}, \ref{MonBeta} and Table \ref{Table:comp} suggest that an analytical estimate will be model specific and parameter dependent.

In \cite{Whittle1955}, Whittle determined that the probability of extinction for a Susceptible-Infected (SI) model was the reciprocal of $\mathcal{R}_0$.  This result was also verified by \cite{Allen2012}.  \cite{Allen2013} showed that $1-\sigma$ and $1-\mathcal{R}_0$ have the same sign, where $\sigma$ is the spectral radius of the matrix of first moments, $\mathbb{M}$.  This implies that efforts that reduce $\mathcal{R}_0$ will also increase the probability of extinction.  For the models studied in this article, one way to reduce $\mathcal{R}_0$ is to decrease the birth rate of susceptible fish, $\beta$.  Unfortunately, this also has the effect of reducing the disease-free equilibrium population size.  Never-the-less, in Figures \ref{MABeta}, \ref{IdBeta} and \ref{MonBeta}, we see that as $\beta$ is decreased, the probability of extinction increases as measured both by branching process approximation and computer simulation.

Metapopulation models are characterized by their patch structure and the rates of migration between patches.  In order to study stochastic metapopulations, it would be useful to study how statistics like probability of extinction vary from patch to patch.  In addition, the probability of partial extinction events, like extinction in one patch, may be useful in measuring the effectiveness of control strategies.  One would expect it to be especially useful when studying the effectiveness of control strategies that are deployed heterogeneously.  Mathematically, the problem of calculating the probability of extinction corresponds to the classical problem of hitting a subspace of the state space of the CTMC from some initial state.  In the case of total extinction events, this relates to hitting the subspace of the state space where all infectious classes are zero.  Taking the two-patch model \eqref{6d} as an example, total disease extinction relates to hitting the subspace $\{S_1,S_2\geq0,I_1=I_2=V_1=V_2=0\}$.  However, for partial extinction events, it relates to hitting a subspace of the state space where some infectious classes are zero, but others are positive.  For example, extinction in patch one of the two-patch model relates to hitting the subspace $\{S_1,S_2\geq0,I_2,V_2>0,I_1=V_1=0\}$.  MTBP's track only the infectious classes and are constructed to calculate the probability of hitting the origin, $\mathbf{0}$.  As such, MTBP approximation is only suited to calculating the probability of total extinction.

\section*{Acknowledgements} This work was conducted with the support from NSF grants DMS-1411853, DMS-1515661 and the Center for Applied Mathematics at University of Florida.  The author would like to thank the referees for their helpful suggestions.

\bibliography{bpapprox}
\bibliographystyle{apalike}

\end{document}